\documentclass[12pt,twoside]{amsart}

\theoremstyle{definition}

\theoremstyle{remark}

\numberwithin{equation}{section}

\date{}

\begin{document}

\centerline{\bf Int. Journal of Math. Analysis, Vol. 4, 2010, no.
15, 713 - 720}

\centerline{}

\centerline{}

\centerline {\Large{\bf Proximinal and \u{C}EBY\u{S}EV Sets }}

\centerline{}

\centerline{\Large{\bf in Normed Linear Spaces}}

\centerline{}

\centerline{\bf {Hadi Haghshenas}}

\centerline{}

\centerline{Department of Mathematics, Birjand University,
Birjand, Iran}

\centerline{h$_{-}$haghshenas60@yahoo.com}

\centerline{}

\centerline{}

{\bf Abstract.} In this paper, we study a part of approximation
theory that present the conditions under which a closed set in a
normed linear space is proximinal or \u{C}eby\u{s}ev.

\centerline{}



\section{Basic Definitions and Preliminaries}
In this section, we collect some definitions which will help us to
describe our results in detail. As the first step, let us fix our
notation. Through this paper, $K$ denotes a non-empty subset of
real normed linear space $(X,\|.\|)$ with the topological dual
space $X^{*}$, $S(X)=\{x\in X;\|x\|=1\}$, $B[x;r]=\{y\in X;
\|y-x\| \leq r\}$ and $B(X)=B[0;1]$.\\For an element $x \in X$, we
define the distance function $d_{K}: X \rightarrow \Bbb{R}$ by
$d_{K}(x)= inf \{\|y-x\| ; y\in K \}$. It is easy to see that the
value of $d_{K}(x)$ is zero if and only if $x$ belongs to
$\overline{K}$, the closure of $K$. The subset $K$ is called
proximinal (resp. \u{C}eby\u{s}ev), if for each $x \in X \setminus
K$, the set of best approximations to $x$ from $K$
$$P_{K}(x)= \{y \in K; \|y-x\|=d_{K}(x)\},$$is nonempty (resp. a singleton). This concept was introduced by S. B. Stechkin
and named after the founder of best approximation theory,
\u{C}eby\u{s}ev.\\It is interesting to know the sufficient
conditions for a subset $K$ of a given normed linear space to be a
proximinal or a \u{C}eby\u{s}ev set, and this is what we want to
consider in this paper.\\ It is not difficult to show that every
proximinal subset $K$ of $X$ is also closed. Now, we state and
prove a sufficient condition for proximinality:\\If $K$ is a
closed subset of a finite-dimensional space $X$, then $K$ is
proximinal. To see this, suppose that $x_{0} \in X \backslash K$
and $r_{0}=d_{K}(x_{0})$. If $r > r_{0}$, then there exists $y \in
K$ such that $\|x_{0}-y\|<r$. Therefore $y \in B[x_{0};r] \bigcap
K$. It follows that $B[x_{0};r] \bigcap K \neq \phi$. If $B_{n}=B[
x_{0};r_{0}+\frac{1}{n}] \bigcap K$, then it is clear that $B_{n}$
is a non-empty compact subset of $X$ and $B_{n+1} \subseteq B_{n}$
for all $n \geq 1$. Hence, there exists $y_{0} \in X$ such that
$y_{0} \in \displaystyle{\ \bigcap_{n=1}^{\infty}B_{n}}$. Now, we
have $\|y_{0}-x_{0}\| \leq r_{0}+\frac{1}{n}$ for all $n \geq 1$.
Since $r_{0}=d_{K}(x_{0})$ we have
$\|y_{0}-x_{0}\|=r_{0}=d_{K}(x_{0})$. Thus $y_{0}$ is a best
approximation for $x_{0}$ and therefore $K$ is a proximinal
set.\\In general, since the functional $e_{x}:
K\rightarrow\Bbb{R}$ with $e_{x}(y)=\|y-x\|$ is continuous, each
compact subset of $X$ is proximinal.\\It is easy to see that in a
reflexive space, every weakly closed set is proximinal.\\
\textbf{Question.   }Is there a closed nonempty subset $K$ of a
reflexive Banach space $X$ with the property that no point outside
$K$ admits a best approximation in $K$? Is this possible in an
equivalent renorm of a Hilbert space? The Lau-Konjagin theorem
(see \cite{2}) states that in a reflexive Banach space $X$, for
every closed set $K$ there is a dense set in $X \setminus K$ which
admits best approximations if and only if the norm has the
Kadec-Klee property. (i.e., for each sequence
$(x_{n})_{n=1}^{\infty}$ in $X$ which converges weakly to $x$ with
$\displaystyle{\lim_{n\rightarrow \infty}\|x_{n}\|= \|x\|}$, we
have $\displaystyle{\lim_{n\rightarrow \infty}\|x_{n}-x\|= 0}$).
\\Every closed convex set in a reflexive space is proximinal
\cite{2}. However, this theorem is not true in the absence of
reflexivity. In fact, this condition is a sufficient one. See the
following example:\\Let $X=l^{1}$. It is known that $l^{1}$ is a
non-reflexive Banach space with dual space $l^{\infty}$. For any
positive integer $n$, let $e_{n} \in l^{1}$ be such that its $n$th
entry is $\frac{(n+1)}{n}$ and all other entries are $0$. Let
$K=\overline{co}\{e_{1},e_{2}, ..., e_{n}, ... \}$. Then $K$ is a
closed convex subset of $l^{1}$ and is not proximinal.\\Another
important notion in this paper is metric projection. The metric
projection mapping has been used in many areas of mathematics such
as the theory of optimization and approximation, and fixed point
theory. It is a set-valued mapping $P_{K}: X \rightarrow K$ which
associates to each $x$ in $X$ the set of all its best
approximations, namely $P_{K}(x)$. The sequence
$(y_{n})_{n=1}^{\infty} \subseteq K$ is called minimizing for $x
\in X \backslash K$ if, $\displaystyle{\lim_{n\rightarrow
\infty}\|x-y_{n}\|= d_{K}(x)}$ and we say that the metric
projection $P_{K}$ is continuous at $x \in X \backslash K$
provided that $\displaystyle{\lim_{n\rightarrow
\infty}y_{n}=y_{0}}$ if, $y_{n} \in P_{K}(x_{n})$ for each $n \in
\Bbb{N}$ and $\displaystyle{\lim_{n\rightarrow \infty}x_{n}=x}$.
It is clear that $P_{K}$ is continuous at $x$ if, every minimizing
sequence for $x \in X \backslash K$ converges \cite{10}. The
continuity properties of $P_{K}$ is a natural object of study in
understanding the nature of some problems in approximation theory.
In the linear cases many results show the connection of the
continuity properties and the geometry of the Banach space (see
\cite{12}). We use this property to prove our main result.\\In
order to give sufficient conditions for a set being proximinal, N.
V. Efimov and S. B. Stechkin introduced the concept of
approximatively compact sets. The set $K$ is said to be
approximatively compact if, for any $x \in X$, each minimizing
sequence $(y_{n})_{n=1}^{\infty} \subseteq K$ for $x$ has a
subsequence converging to an element of $K$. It is proved in
\cite{12} that every approximatively compact set is proximinal.
But a proximinal set need not be approximatively compact
\cite{13}.
\\We say that $K$ is boundedly compact, provided that $K \bigcap
B[0;r]$ is compact in $X$ for every $r\geq 0$. Every boundedly
compact set is approximatively compact although the converse is
false. Thus, every boundedly compact set is proximinal, too. It is
easy to verify that if $K$ is a boundedly compact \u{C}eby\u{s}ev set
in $X$, then the metric projection of $X$ on to $K$ is continuous. Hence,
in a finite dimensional space, every \u{C}eby\u{s}ev set has a continuous metric projection.\\
Let $f:X\rightarrow \Bbb{R}$ be a function and $x,y\in X$. Then
$f$ is said to be Gateaux differentiable at $x$ if, there exists
$A\in X^{*}$ such that $A(y)=\displaystyle{\lim_{t\rightarrow
0}\frac{f(x+ty)-f(x)}{t}}$. In this case $A$ is called the Gateaux
derivative of $f$ and is denoted by $f'(x)$. Also, $A(y)$ is
denoted by $<f'(x),y>$, usually. If the above limit exists
uniformly for each $y\in S(X)$, then $f$ is said to be Fr\'{e}chet
differentiable at $x$ with Fr\'{e}chet derivative $A$. Similarly,
the norm function $\|.\|$ is Gateaux (Fr\'{e}chet) differentiable
at $0 \neq x \in X$ if, the function $f(x)=\|x\|$ is Gateaux
(Fr\'{e}chet) differentiable.\\It is well known that if,
$f:X\rightarrow \Bbb{R}$ is Fr\'{e}chet differentiable at $x \in
X$ then for given $\varepsilon > 0$ there exists
$\delta_{(x,\varepsilon)}> 0$ such that
$\|f(x+y)-f(x)-<f'(x),y>\|\leq\varepsilon\|y\|$, for each $y \in X
$ with $\|y\|< \delta$.\section{Main Results}We start our work
with the following lemma:\\ \textbf{Lemma 1.   }Suppose $K$ is
closed and the distance function  $d_{K}$ is Gateaux
differentiable at $x \in X \backslash K$. Then for every $y \in
P_{K}(x)$ we have $\langle
d'_{K}(x),\frac{x-y}{\|x-y\|}\rangle=1$.\begin{proof}At first,
from Gateaux differentiability of $d_{K}$, the limit
$$\displaystyle{\liminf_{t\rightarrow
0^{+}}\frac{d_{K}(x+tz)-d_{K}(x)}{t}},$$ exists for every $z \in
X$. But for each $t > 0$$$d_{K}(x+t(x-y))-d_{K}(x)\leq t
d_{K}(x).$$Hence, in particular, for
$z=x-y$$$\displaystyle{\liminf_{t\rightarrow
0^{+}}\frac{d_{K}(x+tz)-d_{K}(x)}{t}=d_{K}(x)}.$$Now if
$t'=\frac{t}{d_{K}(x)}$ (notice that $d_{K}(x)>0$ ) then
$$\displaystyle{\liminf_{t'\rightarrow
0^{+}}\frac{d_{K}(x+t'(x-y))-d_{K}(x)}{t'}=d_{K}(x)},$$ and
consequently$$\displaystyle{\liminf_{t\rightarrow
0^{+}}\frac{d_{K}(x+t\frac{x-y}{\|x-y\|})-d_{K}(x)}{t}=1}.$$On the
other hand, since distance functions are Lipschitz (with constant
1) we have$$\displaystyle{\limsup_{t\rightarrow
0^{+}}\frac{d_{K}(x+t\frac{x-y}{\|x-y\|})-d_{K}(x)}{t}\leq1},$$as
required.\end{proof}We say that a non-zero element $x^{*} \in
X^{*}$ strongly exposes $B(X)$ at $x \in S(X)$, provided a
sequence $(z_{n})_{n=1}^{\infty}$ in $B(X)$ converges to $x$
whenever $(\langle x^{*},z_{n}\rangle)_{n=1}^{\infty}$ converges
to $\langle x^{*},x\rangle$.\\The following theorem is the same as
Theorem 2.6 in \cite{10}, but with some manipulation, and plays a
key role in our work:\\\textbf{Theorem 2. }Suppose $K$ is closed
in $X$ and $d_{K}$ is Fr\'{e}chet differentiable at $x \in X
\backslash K$. Moreover $y \in P_{K}(x)$ and $d'_{K}(x)$ strongly
exposes $B(X)$ at $\|x-y\|^{-1}(x-y)$. Then every minimizing
sequence $(y_{n})_{n=1}^{\infty}$ in $K$ for $x$ converges to
$y$.\begin{proof}We can choose a sequence $(a_{n})_{n=1}^{\infty}$
of positive numbers such that $\displaystyle{\lim_{n\rightarrow
\infty}a_{n}=0}$ and
$$a_{n}^{2}>\|x-y_{n}\|-d_{K}(x) \hspace{2 cm}(n \in
\Bbb{N}).$$Hence, if $0 < t < 1$ then for each $n \in \Bbb{N}$
\begin{eqnarray*}
d_{K}(x+t(y_{n}-x))& \leq &\|x+t(y_{n}-x)-y_{n}\| \\& = &(1-t)\|x-y_{n}\|\\
 & < & (1-t)(a_{n}^{2}+d_{K}(x)).
\end{eqnarray*}
Therefore $$d_{K}(x)-d_{K}(x+t(y_{n}-x))\geq
td_{K}(x)-2a_{n}^{2}.$$Fix $\varepsilon > 0$. By Fr\'{e}chet
differentiability of $d_{K}$, there is $\delta > 0$ such that if
$\|y\| < \delta$ then $$|d_{K}(x+y)-d_{K}(x)-\langle
d'_{K}(x),y\rangle|\leq \varepsilon \|y\|\hspace{1 cm}(*).$$Let
$t_{n}=\frac{a_{n}}{\|x-y_{n}\|}$ and $a_{n}< \delta$ for large
$n$. Replacing $y$ by $t_{n}(y_{n}-x)$ in $(*)$ we get
\begin{eqnarray*}
\varepsilon t_{n} \|x-y_{n}\| - \langle
d'_{K}(x),t_{n}(y_{n}-x)\rangle & \geq &
d_{K}(x)-d_{K}(x+t_{n}(y_{n}-x))\\& \geq & t_{n}d_{K}(x)- 2
a_{n}^{2},
\end{eqnarray*}whence$$\langle
d'_{K}(x),t_{n}(x-y_{n})\rangle \geq - \varepsilon a_{n} - 2
a_{n}^{2}+t_{n} d_{K}(x), $$therefore$$\langle
d'_{K}(x),\|x-y_{n}\|^{-1}(x-y_{n})\rangle \geq - \varepsilon - 2
a_{n} + \frac{d_{K}(x)}{\|x-y_{n}\|}.$$Since $\varepsilon > 0$,
$\displaystyle{\lim_{n\rightarrow \infty}a_{n}=0}$,
$\displaystyle{\lim_{n\rightarrow \infty} \|x-y_{n}\|= d_{K}(x)}$,
we will have$$1 \geq \displaystyle{\liminf_{n\rightarrow \infty}
\langle d'_{K}(x),\|x-y_{n}\|^{-1}(x-y_{n}) \rangle}\geq
\displaystyle{\liminf_{n\rightarrow
\infty}\frac{d_{K}(x)}{\|x-y_{n}\|}}=1,$$therefore by the Lemma 1
$$\displaystyle{\lim_{n\rightarrow \infty} \langle
d'_{K}(x),\|x-y_{n}\|^{-1}(x-y_{n})\rangle}=1= \langle
d'_{K}(x),\|x-y\|^{-1}(x-y)\rangle.$$Since $d'_{K}(x)$ strongly
exposes $B(X)$ at $\|x-y\|^{-1}(x-y)$, we deduce that
$$\displaystyle{\lim_{n\rightarrow \infty}
\|x-y_{n}\|^{-1}(x-y_{n})}=\|x-y\|^{-1}(x-y),$$ which yields
$\displaystyle{\lim_{n\rightarrow \infty}y_{n}=y}$.\end{proof}It
is interesting to know that if $K$ is closed in $X$, $x \in X
\backslash K$ and $(y_{n})_{n=1}^{\infty}$ is a minimizing
sequence in $K$ for $x$ with the weak limit $y \in K$, then $y$ is
a best approximation for $x$ in $K$. This is because the norm is a
lower-semi-continuous function with respect to weak topology and
we have$$d_{K}(x)\leq\|x-y\|\leq
\displaystyle{\liminf_{n\rightarrow \infty}\|x-y_{n}\|}\leq
\displaystyle{\lim_{n\rightarrow
\infty}\|x-y_{n}\|}=d_{K}(x).$$\textbf{Theorem 3. }\cite{6} The
dual norm of $X^{*}$ is Fr\'{e}chet differentiable at $x^{*} \in
X^{*}$ if and only if $x^{*}$ strongly exposes
$B(X)$.\\\textbf{Corollary 4. } Let $K$ is closed in $X$, the
distance function $d_{K}$ is Fr\'{e}chet differentiable at $x \in
X \backslash K$ and the dual norm of $X^{*}$ is Fr\'{e}chet
differentiable. Then each minimizing sequence in $K$ for $x$ is
convergent.\begin{proof}Combine Theorem 2 with Theorem
3.\end{proof}\textbf{Corollary 5. } Let $K$ be closed in $X$, $x
\in X \backslash K$ and the distance function $d_{K}$ is
Fr\'{e}chet differentiable at $x$. Also, assume that the dual norm
of $X^{*}$ is
Fr\'{e}chet differentiable. Then the metric projection $P_{K}$ is continuous at $x$.\\

We say that the space $X$ is strictly convex (rotund) if, $x=y$
whenever $\|x\|= \|y\|=\|\frac{x+y}{2}\|=1$ and $X$ is called
uniformly convex if, for sequences $(x_{n})_{n=1}^{\infty},
(y_{n})_{n=1}^{\infty} \subseteq X$ with
$$\displaystyle{\lim_{n\rightarrow
\infty}2\|x_{n}\|^{2}+2\|y_{n}\|^{2}-\|x_{n}+y_{n}\|^{2}}=0,$$we
have $$\displaystyle{\lim_{n\rightarrow
\infty}\|x_{n}-y_{n}\|}=0.$$Obviously, uniformly convex Banach
spaces are strictly convex and also, they are reflexive
(Milman-Pettis).\\\textbf{Remark 6. }It is a well known theorem
that the dual norm of $X^{*}$ is Fr\'{e}chet differentiable if and
only if $X$ is uniformly convex. Therefore, we have the following
corollary.\\\textbf{Corollary 7. }Suppose that $K$ is closed in a
uniformly convex space $X$, $x \in X \backslash K$ and $d_{K}$ is
Fr\'{e}chet differentiable at $x$. Then the metric projection
$P_{K}$ is continuous at $x$.\\We can say also about weakly closed
sets that, each weakly closed set in a uniformly convex Banach
space has continuous metric projection \cite{8}.\\It is proved
that closed convex sets in strictly convex reflexive Banach spaces
(and consequently in uniformly convex Banach spaces) are
\u{C}eby\u{s}ev (see \cite{6}). Can we prove that in some Banach
spaces, a nonempty subset is a \u{C}eby\u{s}ev set if and only if
it is closed and convex? This is an open problem, even in the
special case of infinite-dimensional Hilbert space (see \cite{4}).
In 1934, L. N. H. Bunt proved that each \u{C}eby\u{s}ev set in a
finite-dimensional Hilbert space must be convex. From this result,
we see that in a finite-dimensional Hilbert space, a nonempty
subset is a \u{C}eby\u{s}ev set if and only if it is closed and
convex. In \cite{11}, G. G. Johnson gave an example: there exists
an incomplete inner product space which possesses a non-convex
\u{C}eby\u{s}ev set (M. Jiang completed the proof in 1993). Is
there an infinite-dimensional Hilbert space possessing a
non-convex \u{C}eby\u{s}ev set? As addressed above, it is unknown.
Now, in the last part of the paper, we present a condition under
which a closed subset is \u{C}eby\u{s}ev.
\\It can be seen in \cite{8} that if the dual norm of $X^{*}$ is Fr\'{e}chet
differentiable, then the closed sets in $X$ with continuous
metric projection are \u{C}eby\u{s}ev.\\
Finally, the following is immediate from Corollary 7.
\\\textbf{Corollary 8. }Let $K$ be closed in a uniformly convex
Banach space $X$, $x \in X \backslash K$ and $d_{K}$ is
Fr\'{e}chet differentiable at $x$. Then $K$ is \u{C}eby\u{s}ev in
$X$.\\Corollary 8 is also valid in infinite-dimensional Hilbert
spaces.
\section{Acknowledgments}
The author is indebted to his supervisor professor Amanollah
Assadi for the useful remarks while this work was in progress.
Also, Helpful comments by Mr. H. Hosseini Guive is gratefully
appreciated.

\centerline{}

{\bf Received: September, 2009}
\end{document}